\documentclass[leqno, a4paper]{amsart}

\usepackage[mathscr]{eucal}
\usepackage{amsmath,amsfonts,amssymb}
\usepackage{mathrsfs}
\usepackage[dvips]{graphicx}
\usepackage{bm}
\usepackage{fancyhdr}

\newcommand{\mathsym}[1]{{}}

\theoremstyle{plain}
\newtheorem{theorem}{Theorem}[section]

\newtheorem{proposition}{Proposition}[section]
\newtheorem{lemma}{Lemma}[section]

\theoremstyle{definition}
\newtheorem{definition}{Definition}[section]
\newtheorem{remark}{Remark}[section]

\renewcommand{\footrulewidth}{0pt}

\begin{document}
\title{On the zeros of certain Poincar\'{e} series for $\bm{\Gamma_0^{*} (2)}$ and $\bm{\Gamma_0^{*} (3)}$}
\author{Junichi Shigezumi}
\maketitle \vspace{-0.2in}

\begin{abstract}
We determine the location of all of the zeros of certain Poincar\'{e} series associated with the Fricke groups $\Gamma_0^{*}(2)$ and $\Gamma_0^{*}(3)$ in their fundamental domains by applying and extending the method of F. K. C. Rankin and H. P. F. Swinnerton-Dyer (``{\it On the zeros of Eisenstein series}'', 1970).
\end{abstract}

\thispagestyle{fancy}
\lhead{}\chead{}\rhead{}
\lfoot{{\small {\bfseries 2000 Mathematics Subject Classification:} Primary 11F11.}}
\renewcommand{\footrulewidth}{0.4pt}

\section{Introduction}\label{sec-intro}

F. K. C. Rankin and H. P. F. Swinnerton-Dyer considered the problem of locating the zeros of the Eisenstein series $E_k(z)$ in the standard fundamental domain $\mathbb{F}$ \cite{RSD}. They proved that all of the zeros of $E_k(z)$ in $\mathbb{F}$ are on the unit circle. They also stated towards the end of their study that  ``This method can equally well be applied to Eisenstein series associated with subgroup of the modular group.'' However, it seems unclear how generally this claim holds.

Furthermore, R. A. Rankin considered the same problem for certain Poincar\'{e} series associated with $\text{SL}_2(\mathbb{Z})$ \cite{R}. He also applied the method of F. K. C. Rankin and H. P. F. Swinnerton-Dyer, and proved that all of the zeros of certain Poincar\'{e} series in $\mathbb{F}$ also lie on the unit circle. Also, there are some families of modular forms and fuctions which have similar location of the zeros \cite{Ka}.

Subsequently, T. Miezaki, H. Nozaki, and the present author considered the same problem for Fricke groups $\Gamma_0^{*}(2)$ and $\Gamma_0^{*}(3)$ (See \cite{Kr}, \cite{Q}), which are commensurable with $\text{SL}_2(\mathbb{Z})$. For a fixed prime $p$, we define $\Gamma_0^{*}(p) := \Gamma_0(p) \cup \Gamma_0(p) \: W_p$, where $\Gamma_0(p)$ is a congruence subgroup of $\text{SL}_2(\mathbb{Z})$ and $W_p := \left( \begin{smallmatrix} 0 & -1/\sqrt{p} \\ \sqrt{p} & 0 \end{smallmatrix} \right)$. Let $E_{k, p}^{*}(z)$ be the Eisenstein series associated with $\Gamma_0^{*}(p)$, and the region 
\begin{equation*}
\mathbb{F}^{*}(p) := \left\{|z| \geqslant 1 / \sqrt{p}, \: - 1 / 2 \leqslant Re(z) \leqslant 0\right\}
 \bigcup \left\{|z| > 1 / \sqrt{p}, \: 0 \leqslant Re(z) < 1 / 2 \right\}
\end{equation*}
be a fundamental domain for $\Gamma_0^{*}(p)$ for $p = 2, \; 3$. The authors applied the method of F. K. C. Rankin and H. P. F. Swinnerton-Dyer to the Eisenstein series associated with $\Gamma_0^{*}(2)$ and $\Gamma_0^{*}(3)$. We proved that all of the zeros of $E_{k, p}^{*}(z)$ in $\mathbb{F}^{*}(p)$ lie on the arc $A_p^{*} := \mathbb{F}^{*}(p) \cap \{ |z| = 1 / \sqrt{p} \}$ for $p = 2, 3$ \cite{MNS}.

Let
\begin{equation}
G_{k, \; p}^{*} (z ; R) :=
\frac{1}{2} \sum_{\begin{subarray}{c} (c, d)=1\\ p \mid c\end{subarray}} \frac{R(e^{2 \pi i \gamma z})}{(c z + d)^k} + \frac{1}{2} \sum_{\begin{subarray}{c} (c, d)=1\\ p \mid c\end{subarray}} \frac{R(e^{2 \pi i \gamma (- 1 / (p z))})}{(d \sqrt{p} z - c / \sqrt{p})^k} \label{def-gp}
\end{equation}
be a Poincar\'{e} series associated with $\Gamma_0^{*}(p)$ where $R(t)$ is a suitably chosen rational function of $t$. Here, $\gamma$ is an element of $\Gamma_0(p)$ which satisfies $\gamma = \left(\begin{smallmatrix} a & b \\ c & d \end{smallmatrix}\right)$ for the integers $c, d$ and some integers $a, b$, and is such that
\begin{equation*}
\gamma z = \frac{a z + b}{c z + d}.
\end{equation*}

In the present paper, we consider the same problem for certain Poincar\'{e} series for $\Gamma_0^{*}(2)$ and $\Gamma_0^{*}(3)$. We apply both the method of F. K. C. Rankin and H. P. F. Swinnerton-Dyer, and also the method of R. A. Rankin. The dimension of cusp forms for $\Gamma_0^{*}(p)$ of weight $k$ is denoted by $l$. We prove the following theorems:

\begin{theorem}
Let $k \geqslant 4$ be an even integer, and $m$ be a non-negative integer. Then all of the zeros $(${\it i.e.} $k ((p + 1) / 24) + m$ zeros$)$ of $G_{k, p}^{*}(z ; t^{- m})$ in $\mathbb{F}^{*}(p)$ lie on the arc $A_p^{*}$ for $p = 2, 3$. \label{th-t-m}
\end{theorem}

\begin{theorem}
Let $k \geqslant 4$ be an even integer, and let $m \leqslant l$ be a positive integer. Then $G_{k, p}^{*}(z ; t^m)$ has at least $k ((p + 1) / 24) - m$ zeros on the arc $A_p^{*}$ and at least one zero at $\infty$ for $p = 2, 3$. \label{th-t+m}
\end{theorem}

Note that $G_{k, p}^{*}(z ; t^{- 0}) = E_{k, p}^{*}(z)$.

Furthermore, if a modular form for $\Gamma_0^{*}(p)$ of weight $k$ has $N$ zeros and $P$ poles in $\mathbb{F}^{*}(p)$, then
\begin{equation}
N - P = k ((p + 1) / 24). \quad \text{(See \cite{RSD})}
\end{equation}
In Theorem \ref{th-t-m}, $G_{k, p}^{*}(z ; t^{- m})$ has $m$ poles at $\infty$. On the other hand, in Theorem \ref{th-t+m}, $G_{k, p}^{*}(z ; t^m)$ has no poles and has $k ((p + 1) / 24)$ zeros in total, thus the location of $m - 1$ zeros is unclear.\\

\lhead[\thepage]{\quad}\rhead[\quad]{\thepage}
\chead[{\upshape J. Shigezumi}]{\upshape Zeros of Poincar\'{e} series for Fricke groups}
\lfoot{}
\renewcommand{\footrulewidth}{0pt}

\section{Distribution of the zeros of modular functions}
As is well known, there are two interesting series of modular functions for $\text{SL}_2 (\mathbb{Z})$, for which all of the zeros are on the lower arcs of the fundamental domain $\mathbb{F}$, and with different distributions for the zeros of the two series. They are the Eisenstein series $E_k$ and the Hecke type Faber polynomials $F_m$.

For the former, in the paper of F. K. C. Rankin and H. P. F. Swinnerton-Dyer \cite{RSD}, it was shown that
\begin{equation}
e^{i k \theta / 2} E_k \left(e^{i \theta}\right) = 2 \cos(k \theta / 2) + R_1
\end{equation}
which is real for all $\theta \in [\pi / 2, 2 \pi / 3]$ and also that $|R_1| < 2$ for all $k \geqslant 12$. If $\cos (k \theta / 2)$ is $+1$ or $-1$, then $e^{i k \theta / 2} E_k\left(e^{i \theta}\right)$ is positive or negative, respectively. Thus, the distribution of the zeros of the Eisenstein series resembles a uniform distribution with {\it argument $\theta$} on the lower arc of $\mathbb{F}$.

On the other hand, for the latter, T. Asai, M. Kaneko, and H. Ninomiya proved that
\begin{equation}
F_m(e^{i \theta}) = 2 e^{2 \pi m \sin\theta} \cos(2 \pi m \cos\theta) + R_2
\end{equation}
which is real for all $\theta \in [\pi / 2, 2 \pi / 3]$ and they showed also that $|R_2| < 2 e^{2 \pi m \sin\theta}$ for all $m \geqslant 0$ \cite{AKN}. Thus, the distribution of the zeros resembles a uniform distribution with {\it real part $\sin\theta$} on the lower arc of $\mathbb{F}$.

Now, we consider the Poincar\'{e} series (cf. Eq. (\ref{def-gp})):
\begin{equation}
G_k (z ; t^{-m}) := 
\frac{1}{2} \sum_{(c, d)=1} \frac{R(e^{2 \pi i \gamma z})}{(c z + d)^k}.
\end{equation}
In the paper of R. A. Rankin \cite{R}, it was shown that
\begin{equation}
e^{i k \theta / 2} G_k \left( e^{i \theta} ; t^{-m} \right) = 2 e^{2 \pi m \sin\theta} \cos(k \theta / 2 - 2 \pi m \cos\theta) + R_3
\end{equation}
which is real for all $\theta \in [\pi / 2, 2 \pi / 3]$ and also that $|R_3| < 2 e^{2 \pi m \sin\theta}$ for all $k \geqslant 4$ and $m \geqslant 0$. Then, Rankin proved that all of the zeros of $G_k (z; \: t^{-m})$ are on the lower arc of $\mathbb{F}$.

Note that $G_k (z ; t^0) = E_k (z)$ and $G_0 (z ; t^{-m}) = F_m (z)$. Furthermore, if $k$ is large enough compared with $m$, then the distribution of zeros of $G_k (z ; t^{-m})$ resembles that of $E_k (z)$. On the other hand, if $m$ is large enough compared with $k$, then the distribution resembles that of $F_m (z)$. Thus, the Poincar\'{e} series $G_k (z ; t^{-m})$ ``fill the space of two modular functions discretely''. For example, we consider the following sequence of modular forms:
\begin{equation}
G_{12 l} (z ; t^0), \: G_{12 (l - 1)} (z ; t^{- 1}), \: G_{12 (l - 2)} (z ; t^{- 2}), \: \ldots, \: G_0 (z ; t^{-l}).
\end{equation}
The number of zeros of each function is equal to $l$, and all of the zeros are on the lower arc of $\mathbb{F}$.\\

For the Fricke groups $\Gamma_0^{*}(2)$ and $\Gamma_0^{*}(3)$, the distributions of zeros are similar. In fact, we have the following relations:
\begin{align*}
e^{i k \theta / 2} E_{k, p}^{*} \left( e^{i \theta} / \sqrt{p} \right) &= 2 \cos(k \theta / 2) + R_1,\\
F_{m, p}(e^{i \theta} / \sqrt{p}) &= 2 e^{(2 / \sqrt{p}) \pi m \sin\theta} \cos((2 / \sqrt{p}) \pi m \cos\theta) + R_2,\\
e^{i k \theta / 2} G_{k, p}^{*} \left( e^{i \theta} / \sqrt{p} ; t^{-m} \right) &= 2 e^{(2 / \sqrt{p}) \pi m \sin\theta} \cos(k \theta / 2 - (2 / \sqrt{p}) \pi m \cos\theta) + R_3,
\end{align*}
and it has been shown that all the zeros of the functions are on the lower arcs of $\mathbb{F}^{*}(p)$ for the Eisenstein series $E_{k, p}^{*}$ (cf. \cite{MNS}), the Hecke type Faber polynomials $F_{k, 2}$ (cf. \cite{BKM}), and the Poincar\'{e} series $G_{k, p}^{*}$ (cf. the present paper).

\begin{remark}
W. Duke and P. Jenkins considered certain weakly holomorphic modular forms $f_{k, m}$ for $\text{SL}_2 (\mathbb{Z})$. If we assume that $k \geqslant 4$ and $m \geqslant 0$ then $f_{k, m}$ is a modular function of weight $k$ and with $m$ poles at $\infty$. They proved that all the zeros of the function are on the lower arc of $\mathbb{F}$. Note that we have the following relation:
\begin{equation*}
e^{i k \theta / 2} f_{k, m} \left( e^{i \theta} \right) = 2 e^{2 \pi m \sin\theta} \cos(k \theta / 2 - 2 \pi m \cos\theta) + R.
\end{equation*}
Thus, the distribution of the zeros is similar to that of the Poincar\'{e} series $G_k (z ; t^{-m})$.

However, unfortunately, for $\Gamma_0^{*}(2)$ and $\Gamma_0^{*}(3)$, not all of the zeros of the function $f_{k, m}$ are on the lower arc of $\mathbb{F}^{*}(p)$. H. Tokitsu observed it by numerical calculation \cite{T}, and he proved that if $m$ is large enough then all of the zeros of $f_{k, m}$ are on the lower arc.
\end{remark}\quad

\section{General Theory}
Let $\theta_{0, 2} := 3 \pi / 4$, $\theta_{0, 3} := 5 \pi / 6$, and $\theta_{1, p} := \pi / 2$ for $p = 2, 3$, and let $\rho_p := e^{i \theta_{0, p}} / \sqrt{p}$. We denote by $v_q(f)$ the order of a modular function $f$ at a point $q$. We have the following propositions: (See \cite{MNS})

\begin{proposition}
Let $k \geqslant 4$ be an even integer. For every function $f$ which is a modular form for $\Gamma_0^{*}(p)$, we have
\begin{equation}
\begin{split}
v_{i / \sqrt{p}}(f) \geqslant s_k &\quad(s_k=0, 1 \; \text{such that} \; 2 s_k \equiv k \pmod{4}),\\
v_{\rho_p}(f) \geqslant t_k &\quad(0 \leqslant t_k < 12 / (p + 1) \; \text{such that} \; - 2 t_k \equiv k \pmod{24 / (p + 1)}).
\end{split}
\end{equation} \label{prop-bd_ord_p}
\end{proposition}

Furthermore, let $l$ be the dimension of cusp forms for $\Gamma_0^{*}(p)$ of weight $k$; then we have $l = \left\lfloor k ((p + 1) / 24) - t / 4 \right\rfloor$, where $t = 0 \text{ or } 2$, such that $t \equiv k \pmod{4}$.

Following the methods in \cite{RSD} and \cite{R}, we define
\begin{equation}
F_{k, \; p}^{*} (\theta ; R) := e^{i k \theta / 2} G_{k, \; p}^{*} \left( e^{i \theta} / \sqrt{p} \, ; R \right).
\end{equation}
Then, we have
\begin{equation*}
F_{k, \; p}^{*} (\theta ; R) =
\frac{1}{2} \sum_{\begin{subarray}{c} (c, d)=1\\ p \mid c\end{subarray}} \frac{R \left( e^{2 \pi i \gamma (e^{i \theta} / \sqrt{p})} \right)}{\left( (c / \sqrt{p}) e^{i \theta / 2} + d e^{- i \theta / 2} \right)^k} + \frac{1}{2} \sum_{\begin{subarray}{c} (c, d)=1\\ p \mid c\end{subarray}} \frac{R \left( e^{2 \pi i \gamma (- (e^{- i \theta} / \sqrt{p}))} \right)}{\left( (c / \sqrt{p}) e^{- i \theta / 2} - d e^{i \theta / 2} \right)^k},
\end{equation*}
which is real for every $\theta \in [\theta_{1, p}, \theta_{0, p}]$.

Furthermore, we can write
\begin{equation}
F_{k, \; p}^{*} (\theta ; R) = 2 \, Re \, g_{k, \; p}^{*} (\theta ; R) + F_{k, \; p}^{* '} (\theta ; R),
\end{equation}
where $F_{k, \; p}^{* '}$ consists of all of the terms of the series $F_{k, \; p}^{*}$ which satisfy $c^2 + d^2 \geqslant 2$, and where
\begin{equation*}
g_{k, \; p}^{*} (\theta ; R) = e^{i k \theta / 2} R \left( e^{2 \pi i (e^{i \theta} / \sqrt{p})} \right).
\end{equation*}

Let $\gamma_p$ be the locus of $e^{2 \pi i (e^{i \theta} / \sqrt{p})}$ as $\theta$ increases from $2 \theta_{1, p} - \theta_{0, p}$ to $\theta_{0, p}$, and let $r_{j, p} := |e^{2 \pi i (e^{i \theta_{j, p}} / \sqrt{p})}|$ for $j = 0, 1$. Then, the curve $\gamma_p$ begins at $-r_{0, p}$, passes through $r_{1, p}$, and returns to $-r_{0, p}$ following a clockwise rotation. Assume that the function $R$ has no zero or pole on $\gamma_p$, and that it has $N_{\gamma_p}$ zeros and $P_{\gamma_p}$ poles in the region $\mathbb{D}_{\gamma_p}$, the interior of $\gamma_p$. Then, by the Argument Principle, we have
\begin{equation}
n_0 - n_1 =P_{\gamma_p} - N_{\gamma_p},
\end{equation}
where $n_0$ and $n_1$ are integers such that
\begin{equation*}
\arg g_{k, \; p}^{*} (\theta_{1, p} ; R) = \pi (n_1 + k / 4), \quad \arg g_{k, \; p}^{*} (\theta_{0, p} ; R) = \pi (n_0 + k (\theta_{0, p} / (2 \pi))).
\end{equation*}

Suppose that the Poincar\'{e} series $G_{k, p}^{*}(z ; R)$ has $N_R$ zeros and $P_R$ poles in $\mathbb{F}^{*}(p)$. Then
\begin{equation}
N_R - P_R = k ((p + 1) / 24).
\end{equation}

\begin{definition}[Property $P_{k, p}$]
We shall say that the function $R$ has Property $P_{k, p}$ if $(i)$ $R$ is a real rational function, $(ii)$ all of the poles of $R$ lie in $\mathbb{D}_{\gamma_p}$, $(iii)$ $l \geqslant N_{\gamma_p} - P_{\gamma_p}$, and $(iv)$
\begin{equation}
|F_{k, \; p}^{* '} (\theta ; R)| < 2 \left| R \left( e^{2 \pi i (e^{i \theta} / \sqrt{p})} \right) \right| \label{cond-pkp0}
\end{equation}
for every $\theta \in [\theta_{1, p}, \theta_{0, p}]$.
\end{definition}

Here, condition $(iv)$ is equivalent to condition {\it $(iv)'$: the inequality $(\ref{cond-pkp0})$ holds for every point $\theta \in [\theta_{1, p}, \theta_{0, p}]$ which satisfies}
\begin{equation}
\arg g_{k, \; p}^{*} (\theta ; R) \equiv 0 \pmod{\pi}.
\end{equation}
We call such a point an {\it integer point of $g_{k, \; p}^{*}$}.

We then have the following theorem: (See Theorem 1 and 2 of \cite{R})

\begin{proposition}
Suppose that the function $R$ has Property $P_{k, p}$. Then the Poincar\'{e} series $G_{k, p}^{*}(z ; R)$ has at least $N_R - N_{\gamma_p}$ zeros on the arc $A_p^{*}$ for $p = 2, 3$. In particular, if $R$ does not vanish in $\mathbb{D}_{\gamma_p}$, then all of the zeros of $G_{k, p}^{*}(z ; R)$ lie on $A_p^{*}$.\label{prop-pro-pkp}
\end{proposition}

Note that we can prove the above theorem with condition $(iv)'$ replacing condition $(iv)$.

We would like to put condition (\ref{cond-pkp0}) of Property $P_{k, p}$ into another form. We define the following bounds:
\begin{gather}
M_{R, p} := \sup \left\{ \left| R \left( e^{2 \pi i \gamma (e^{i \theta} / \sqrt{p})} \right) \right| \; ; \; \theta_{1, p} \leqslant \theta \leqslant \theta_{0, p} \; ,\; \gamma \in \Gamma_0(p) \setminus \Gamma_{\infty} \right\},\\
\alpha_{k, p} := \sup \left\{ \sum_{\begin{subarray}{c} (c, d)=1\\ p \mid c, \; c \ne 0 \end{subarray}} \left( c^2 / p + d^2 + (2 / \sqrt{p}) c d \cos\theta \right)^{- k / 2} \; ; \; \theta_{1, p} \leqslant \theta \leqslant \theta_{0, p} \right\}.
\end{gather}
Then the condition
\begin{equation}
M_{R, p} \alpha_{k, p} < 2 \left| R \left( e^{2 \pi i (e^{i \theta} / \sqrt{p})} \right) \right| \label{cond-pkp1}
\end{equation}
is  sufficient to imply condition (\ref{cond-pkp0}). Then, we have only to prove the above inequality instead of condition (\ref{cond-pkp0}).

This idea is due to R. A. Rankin \cite{R}.

However, it is difficult to apply this method to the cases of $\Gamma_0^{*}(2)$ and $\Gamma_0^{*}(3)$. For example, if $R(t) = t^{- m}$ for a positive integer $m$, then we have $M_{R, p} = r_{0, p}^{- m}$ and $\left| R \left( e^{2 \pi i (e^{i \theta} / \sqrt{p})} \right) \right| \geqslant r_{0, p}^{- m}$. Moreover, we have $\alpha_{k, 2} > 2$ and $\alpha_{k, 3} > 4$. Thus, we are unable to prove Theorem \ref{th-t-m} in this way, nor Theorem \ref{th-t+m}. We must consider a certain extension of this method, observing some terms of the series in detail.\\

\section{Applications}

The point in the previous section is that there exist some pairs $(c, d)$ such that $c^2 / p + d^2 + (2 / \sqrt{p}) c d \cos\theta_{0, p} = 1$. For the case $p = 2$, the pair satisfying $c^2 / 2 + d^2 + (2 / \sqrt{2}) c d \cos\theta_{0, 2} = 1$ is $(c, d) = \pm(2, 1)$, then $\gamma = \pm S_2 = \pm \left(\begin{smallmatrix} 1 & 0 \\ 2 & 1 \end{smallmatrix}\right)$, which is shown in the term of the Poincar\'{e} series. For the other case $p = 3$, such pairs are given by $(c, d) = \pm(3, 1)$ [in which case $\gamma = \pm S_3 = \pm \left(\begin{smallmatrix} 1 & 0 \\ 3 & 1 \end{smallmatrix}\right)$] and $(c, d) = \pm(3, 2)$ [in which case $\gamma = \mp S_3^{-1} T = \pm \left(\begin{smallmatrix} -1 & -1 \\ 3 & 2 \end{smallmatrix}\right)$]. Note that ``$c^2 / p + d^2 + (2 / \sqrt{p}) c d \cos\theta < 1 \Leftrightarrow \theta < \theta_{0, p}$ for the above pairs $(c, d)$.''

Now, we can write
\begin{gather}
F_{k, \; 2}^{*} (\theta ; R) = 2 \, Re \, g_{k, \; 2}^{*} (\theta ; R) + 2 \, Re \, h_{k, \; 2}^{*} (\theta ; R) + F_{k, \; 2}^{* ''} (\theta ; R),\\
F_{k, \; 3}^{*} (\theta ; R) = 2 \, Re \, g_{k, \; 3}^{*} (\theta ; R) + 2 \, Re \, h_{k, \; 3}^{*} (\theta ; R) + 2 \, Re \, h_{k, \; 3}^{* '} (\theta ; R) + F_{k, \; 3}^{* ''} (\theta ; R)
\end{gather}
where $F_{k, \; p}^{* ''}$ consists of all terms of the series $F_{k, \; p}^{* ''}$ which satisfy $c^2 + d^2 \geqslant 2$ and are not equal to the above pairs, and where
\begin{equation*}
h_{k, \; p}^{*} (\theta ; R) = \frac{R \left( e^{2 \pi i S_p (e^{i \theta} / \sqrt{p})} \right)}{\left( \sqrt{p} e^{i \theta / 2} + e^{- i \theta / 2} \right)^k}, \quad
 h_{k, \; 3}^{* '} (\theta ; R) = \frac{R \left( e^{2 \pi i (- S_3^{-1} T) (e^{i \theta} / \sqrt{3})} \right)}{\left( \sqrt{3} e^{i \theta / 2} + 2 e^{- i \theta / 2} \right)^k}.
\end{equation*}

Instead of $M_{R, p}$, we define
\begin{align}
{M_{R, 2}}' &:= \sup \left\{ \left| R \left( e^{2 \pi i \gamma (e^{i \theta} / \sqrt{2})} \right) \right| \; ; \; \theta_{1, 2} \leqslant \theta \leqslant \theta_{0, 2} \; ,\; \gamma \in \Gamma_0(2) \setminus (\Gamma_{\infty} \cup \Gamma_{\infty} S_2) \right\},\\
{M_{R, 3}}' &:= \sup \left\{ \left| R \left( e^{2 \pi i \gamma (e^{i \theta} / \sqrt{3})} \right) \right| \; ; \; \theta_{1, 3} \leqslant \theta \leqslant \theta_{0, 3} \; ,\; \gamma \in \Gamma_0(3) \setminus (\Gamma_{\infty} \cup \Gamma_{\infty} S_3 \cup \Gamma_{\infty}S_3^{-1} T) \right\}.
\end{align}
Moreover, since
\begin{equation*}
Im \; \gamma (e^{i \theta} / \sqrt{p}) = (1 / \sqrt{p}) \sin\theta / (c^2 p + d^2 + 2 \sqrt{p} c d \cos\theta),
\end{equation*}
 $2 c^2 + d^2 + 2 \sqrt{2} c d \cos\theta \geqslant \sqrt{2} \sin\theta$, and $3 c^2 + d^2 + 2 \sqrt{3} c d \cos\theta \geqslant 2 \sin\theta$, we have $| e^{2 \pi i \gamma (e^{i \theta} / \sqrt{p})} | \geqslant r_{0, p}$. Then, we have
\begin{align*}
{M_{R, 2}}' &= \sup \left\{ |R(t)| \; ; \; e^{- \sqrt{2} \pi / 3} \leqslant |t| \leqslant 1 \right\},\\
{M_{R, 3}}' &= \sup \left\{ |R(t)| \; ; \; e^{- \pi / (2 \sqrt{3})} \leqslant |t| \leqslant 1 \right\}.
\end{align*}
 
On the other hand, for $\alpha_{k, p}$, since $\left( c^2 / p + d^2 + (2 / \sqrt{p}) c d \cos\theta \right)^{- k / 2} + \left( c^2 / p + d^2 - (2 / \sqrt{p}) c d \cos\theta \right)^{- k / 2}$ is monotonically increasing in $\theta$, we have
\begin{align*}
\alpha_{k, p} &= \sum_{\begin{subarray}{c} (c, d)=1\\ p \mid c, \; c \ne 0 \end{subarray}} ( c^2 / p + d^2 + c d )^{- k / 2} \; = \sum_{\begin{subarray}{c} (c, d)=1\\ p \nmid d, \; c \ne 0 \end{subarray}} ( c^2 p + c d p + d^2 )^{- k / 2}\\
 &\leqslant \sum_{(c, d)=1} ( c^2 p + c d p + d^2 )^{- k / 2} - 2.
\end{align*}
In addition, since we have $2 c^2 + 2 c d + d^2 = c^2 + (c + d)^2$ and $3 c^2 + 3 c d + d^2 = c^2 + c (c + d) + (c + d)^2$, we can regard these series as the Epstein zeta-function.  We therefore have the following bounds:
\begin{equation*}
\alpha_{k, 2} \leqslant \frac{4 \, \zeta (k / 2) Z_4(k / 2)}{\zeta (k)} - 2, \quad
 \alpha_{k, 3} \leqslant \frac{6 \, \zeta (k / 2) Z_3(k / 2)}{\zeta (k)} - 2,
\end{equation*}
where $\zeta$ is the Riemann zeta-function and $Z_3$ and $Z_4$ are the Dirichlet $L$-series
\begin{align*}
Z_3(s) &:= 1 - 2^{- s} + 4^{- s} - 5^{- s} + 7^{- s} - 8^{- s} + \cdots,\\
Z_4(s) &:= 1 - 3^{- s} + 5^{- s} - 7^{- s} + 9^{- s} - 11^{- s} + \cdots.
\end{align*}
Furthermore, if $2 c^2 + 2 c d + d^2 = 2$, then $(c, d) = \pm (1, 0), \pm(1, -2)$. These pairs satisfy $(c, d) = 1$, but they do not satisfy $2 \nmid d$. In addition, by using the bounds
\begin{align*}
\zeta(x) &\leqslant 1 + 2^{-x} + 3^{-x} + 4^{-x} + 5^{-x} + 5^{1 - x} / (x - 1),\\
Z_4(x) &\leqslant 1 - 3^{-x} + 5^{-x},\qquad
 \{\zeta(2 x)\}^{-1} \leqslant 1 - 2^{- 2 x},
\end{align*}
we have
\begin{equation*}
\alpha_{k, 2} - 2 - 4 \times 2^{- k / 2} \leqslant 5^{- k / 2} (5 + 10 / (k - 2)) \quad \text{for} \; k \geqslant 4.
\end{equation*}
Then, we define
\begin{equation}
\delta_{k, 2} := 5^{- k / 2} (5 + 10 / (k - 2)).
\end{equation}
Similarly, we define
\begin{equation}
\delta_{k, 3} := \begin{cases}
&7^{- k / 2} (14 + 14 / (k - 2)) \quad \text{for} \; k \geqslant 4\\
&7^{- k / 2} (13 + 14 / (k - 2)) \quad \text{for} \; k \geqslant 12
\end{cases}.
\end{equation}
Finally, we have
\begin{equation*}
|F_{k, \; p}^{* ''}| < {M_{R, p}}' \delta_{k, p}.
\end{equation*}\quad

Now, we define the following condition instead of ``Property $P_{k, p}$'':
\begin{definition}[Property ${P_{k, p}}'$]
We shall say that the function $R$ has Property ${P_{k, p}}'$ if $(i)$ $R$ is a real rational function, $(ii)$ all of the poles of $R$ lie in $\mathbb{D}_{\gamma_p}$, $(iii)$ $l \geqslant N_{\gamma_p} - P_{\gamma_p}$, and
\begin{trivlist}
\item{$(iv-i)$} We have
\begin{gather}
| Re \, h_{k, \; 2}^{*} (\theta ; R)| + (1/2) {M_{R, 2}}' \delta_{k, 2} < \left| R \left( e^{2 \pi i (e^{i \theta} / \sqrt{2})} \right) \right|, \label{cond-pkpp12}\\
| Re \, h_{k, \; 3}^{*} (\theta ; R)| + | Re \, h_{k, \; 3}^{* '} (\theta ; R)| + (1/2) {M_{R, 3}}' \delta_{k, 3} < \left| R \left( e^{2 \pi i (e^{i \theta} / \sqrt{3})} \right) \right| \label{cond-pkpp13}
\end{gather}
for every integer point $\theta \in [\theta_{1, p}, \theta_{0, p})$ for $p = 2, 3$, respectively.
\item{$(iv-ii)$} We have
\begin{equation}
Sign(F_{k, \; p}^{*} (\theta_{0, p} ; R)) = Sign(Re \, g_{k, \; p}^{*} (\theta_{0, p} ; R))
\end{equation}
if $\theta_{0, p}$ is an integer point. \label{cond-pkpp2}
\end{trivlist}
\end{definition}

To prove Theorem \ref{th-t-m} and \ref{th-t+m}, we consider the following theorem instead of Proposition \ref{prop-pro-pkp}, where the point is to use Property ${P_{k, p}}'$.

\begin{theorem}
Suppose that the function $R$ has Property ${P_{k, p}}'$. Then the Poincar\'{e} series $G_{k, p}^{*}(z ; R)$ has at least $N_R - N_{\gamma_p}$ zeros on the arc $A_p^{*}$ for $p = 2, 3$. In particular, if $R$ does not vanish in $\mathbb{D}_{\gamma_p}$, then all of the zeros of $G_{k, p}^{*}(z ; R)$ lie on $A_p^{*}$.\label{th-pro-pkp}
\end{theorem}\quad

\section{Proof of Theorem \ref{th-t-m}}

\subsection{Preliminaries}
Let $k \geqslant 4$ and $m > 0$ be integers, and let $R(t) = t^{-m}$. Then, it is clear that $R$ satisfies the conditions $(i)$, $(ii)$, and $(iii)$ of Property ${P_{k, p}}'$. Furthermore, we have $P_{\gamma_p} = m$, $N_{\gamma_p} = 0$, and $N_R - N_{\gamma_p} = k ((p + 1) / 24) + m$.

To prove that $R$ satisfies condition $(iv-i)$ of Property ${P_{k, p}}'$, it is sufficient to prove the inequalities (\ref{cond-pkpp12}) and (\ref{cond-pkpp13}) for every $\theta \in [\theta_{1, p}, \theta_{0, p} - x]$ for certain $x$ such that every integer point is included in the interval. The first step is to consider how small $x$ should be.

When $p = 2$, and when $k \equiv 4 \pmod{8}$, then $k = 8 l + 4$ and
\begin{equation*}
\arg g_{k, \; 2}^{*} (\theta_{0, 2} ; R) = ((3 / 8) k + m) \pi = (3 l + m + 1 + 1 / 2) \pi.
\end{equation*}
Thus, the last integer point is the point $\theta_{0, p} - x$ such that
\begin{equation*}
\arg g_{k, \; 2}^{*} (\theta_{0, 2} - x ; R) = (3 l + m + 1 + 1 / 2) \pi - \pi / 2.
\end{equation*}
In addition, we have
\begin{equation*}
\arg g_{k, \; 2}^{*} (\theta_{0, 2} - x ; R) > \arg g_{k, \; 2}^{*} (\theta_{0, 2} ; R) - ((k + 8 m) / 2) x.
\end{equation*}
Thus, $x \leqslant \pi / (k + 8 m)$ is sufficient. Similarly, when $k \equiv 6, 0, \, \text{and} \,  2 \pmod{8}$, we have $x \leqslant \pi / (2 (k + 8 m)), 2 \pi / (k + 8 m)$, and $3 \pi / (2 (k + 8 m))$, respectively. So for this case, the bound $x \leqslant \pi / (2 (k + 8 m))$ is sufficient.

Similarly, when $p = 3$, and when $k \equiv 4, 6, 8, 10, 0, \, \text{and} \,  2 \pmod{12}$, we have $x \leqslant t \pi / (k + 6 m)$, where $t = 4 / 3, 1, 2 / 3, 1 / 3, 2, \, \text{and} \, 5 / 3$, respectively. In conclusion, the bound $x \leqslant \pi / (3 (k + 6 m))$ is sufficient. Note that $x \leqslant 2 \pi / (3 (k + 6 m))$ is sufficient if $k \not\equiv 10 \pmod{12}$.

Furthermore, we have  $\left| R \left( e^{2 \pi i (e^{i \theta} / \sqrt{p})} \right) \right| = e^{ (2 / \sqrt{p}) \pi m \sin\theta}$; we define
\begin{gather*}
G := e^{- \sqrt{2} \pi m \sin\theta} |h_{k, \; 2}^{*} (\theta ; R)|,\\
G_1 := e^{- (2 / \sqrt{3}) \pi m \sin\theta} |h_{k, \; 3}^{*} (\theta ; R)|, \quad
G_2 := e^{- (2 / \sqrt{3}) \pi m \sin\theta} |h_{k, \; 3}^{* '} (\theta ; R)|.
\end{gather*}

Then, we will show the following lemma in the following sections:
\begin{lemma}\label{lem-t-m}\quad
\begin{trivlist}
\item{$(i)$} We have
\begin{gather}
G + (1 / 2) \, {M_{R, 2}}' \, \delta_{k, 2} \, e^{- \pi m} < 1 \quad \text{for every} \; \theta \in [\theta_{1, 2}, \theta_{0, 2} - \pi / (2 (k + 8 m))], \label{cond-pkp-2} \\
G_1 + G_2 + (1 / 2) \, {M_{R, 3}}' \, \delta_{k, 3} \, e^{- \pi m} < 1 \quad \text{for every} \; \theta \in [\theta_{1, 3}, \theta_{0, 3} - x_0], \label{cond-pkp-3}
\end{gather}
where $x_0 = \pi / (3 (k + 6 m))$ when $k \equiv 10 \pmod{12}$ and $x_0 = 2 \pi / (3 (k + 6 m))$ when $k \not\equiv 10 \pmod{12}$.

\item{$(ii)$} We have
\begin{equation}
Sign(F_{k, \; p}^{*} (\theta_{0, p} ; R)) = Sign(Re \, g_{k, \; p}^{*} (\theta_{0, p} ; R))
\end{equation}
if $k \equiv 0 \pmod{8}$ for the case of $p = 2$, and if $k \equiv 0 \pmod{12}$ for the case of $p = 3$.
\end{trivlist}
\end{lemma}

When we have proved the lemma above, then we can show that the function $R$ satisfies Property ${P_{k, p}}'$. Thus, we can prove Theorem \ref{th-t-m}.\\

\subsection{The case $p = 2$}
We have ${M_{R, 2}}' = e^{(\sqrt{2} / 3) \pi m}$ and $(1 / 2) \, {M_{R, 2}}' \, \delta_{k, 2} \, e^{- \pi m} \leqslant 0.038003...$ for $k \geqslant 4$ and $m \geqslant 1$.

Furthermore,
\begin{equation*}
G = \text{Exp}[- 2 \pi m \sin\theta (\sqrt{2} + 2 \cos\theta) / (3 + 2 \sqrt{2} \cos\theta)] \; / \; (3 + 2 \sqrt{2} \cos\theta)^{k / 2}.
\end{equation*}
For $\theta \in [\theta_{1, 2}, \theta_{0, 2} - \pi / (2 (k + 8 m))]$, we have
\begin{align*}
&\text{Exp}[- 2 \pi m \sin\theta (\sqrt{2} + 2 \cos\theta) / (3 + 2 \sqrt{2} \cos\theta)]\\
 &\leqslant \text{Exp}[- (20/13) \pi m \sin(\pi / (2 (k + 8 m)))]
 \leqslant \text{Exp}[- (99/130) \pi^2 / ((k / m) + 8)],
\end{align*}
which is monotonically increasing in $k / m$. We also have
\begin{equation*}
3 + 2 \sqrt{2} \cos(\theta_{0, 2} - \pi / (2 (k + 8 m))) \geqslant 1 + \pi (1 / (1 + 8 m / k)) (1 / k).
\end{equation*}

Put $s := k / m$. For $0 \leqslant s \leqslant 100$, we have
\begin{equation*}
G \leqslant \text{Exp}[- (99/130) \pi^2 / (s + 8)] \leqslant 0.93277...
\end{equation*}
On the other hand, if $100 \leqslant s$, we have $1 + \pi (1 / (1 + 8 m / k)) (1 / k) \geqslant 2^{k / 2}$ and
\begin{equation*}
G \leqslant (1 / 2) \text{Exp}[- (99/130) \pi^2 / (s + 8)] \leqslant 1 / 2.
\end{equation*}

Then, $(i)$ of Lemma \ref{lem-t-m} follows.

On the other hand, when $k \equiv 0 \pmod{8}$, we can write $k = 8 l$, and we have
\begin{gather*}
Re \, g_{8 l, \; 2}^{*} (\theta_{0, 2} ; R) = e^{\pi m} \cos((3 l + m) \pi),\\
Re \, h_{8 l, \; 2}^{*} (\theta_{0, 2} ; R) = e^{\pi m} \cos((l - m) \pi).
\end{gather*}
Note that the signs of the two terms above are in agreement. Furthermore, we have $|F_{8 l, \; 2}^{* ''} (\theta_{0, 2} ; R)| \ll |Re \, g_{8 l, \; 2}^{*} (\theta_{0, 2} ; R)|$.

Thus, we can show Lemma \ref{lem-t-m}.\\

\subsection{The case $p = 3$}
We have ${M_{R, 3}}' = e^{(1 / (2 \sqrt{3})) \pi m}$. Let $D := (1 / 2) \, {M_{R, 3}}' \, \delta_{k, 3} \, e^{- \pi m}$.

Firstly, we have
\begin{align*}
G_1 &= \text{Exp}[- 2 \pi m \sin\theta (\sqrt{3} + 2 \cos\theta) / (4 + 2 \sqrt{3} \cos\theta)] \; / \; (4 + 2 \sqrt{3} \cos\theta)^{k / 2},\\
G_2 &= \text{Exp}[- 4 \pi m \sin\theta (\sqrt{3} + 2 \cos\theta) / (7 + 4 \sqrt{3} \cos\theta)] \; / \; (7 + 4 \sqrt{3} \cos\theta)^{k / 2}.
\end{align*}

For $k \geqslant 4$, we have $D = 0.061157...$.

When $k = 4$, we have
\begin{align*}
G_1 + G_2 \leqslant \text{Exp}[- (5/6) \pi m \sin(4 \pi / (3 (4 + 6 m)))] + \text{Exp}[- (10/7) \pi m \sin(4 \pi / (3 (4 + 6 m)))],
\end{align*}
where the right-hand side is monotonically decreasing in $m$ and is less than $1 - D$ when $m = 1$.

When $6 \leqslant k \leqslant 32$ and $k \not\equiv 10 \pmod{12}$ (i.e. $k \ne 10, 22$), we have  
\begin{align*}
G_1 + G_2 &\leqslant \text{Exp}[- (5/3) \pi m \sin(\theta_{0, 3} - 2 \pi / (3 (k + 6 m))) (\sqrt{3} + 2 \cos(\theta_{0, 3} - 2 \pi / (3 (k + 6 m))))]\\
 &\quad + \text{Exp}[- (20/7) \pi m \sin(\theta_{0, 3} - 2 \pi / (3 (k + 6 m))) (\sqrt{3} + 2 \cos(\theta_{0, 3} - 2 \pi / (3 (k + 6 m))))]\\
 &=: f_1(k, m),
\end{align*}
where $f_1(k, m)$ is monotonically decreasing in $m$. Moreover, we have $G_1 + G_2 + D < 1$ for $m \leqslant 6$ and $f_1(k, 7) + D < 1$ by numerical calculation.

For $k \geqslant 10$, we have $D = 0.000078006...$.

When $k =10$, we have $G_1 + G_2 + D < 1$ for $m \leqslant 12$ and
\begin{align*}
G_1 + G_2 &\leqslant \text{Exp}[- (1000/511) \pi m \sin(\theta_{0, 3} - \pi / (3 (10 + 6 m))) (\sqrt{3} + 2 \cos(\theta_{0, 3} - 2 \pi / (3 (10 + 6 m))))]\\
 &\quad + \text{Exp}[- (1000/261) \pi m \sin(\theta_{0, 3} - \pi / (3 (10 + 6 m))) (\sqrt{3} + 2 \cos(\theta_{0, 3} - 2 \pi / (3 (10 + 6 m))))],
\end{align*}
where the right-hand side is monotonically decreasing in $m$ and is less than $1 - D$ when $m = 13$.

When $k = 22$, we have $G_1 + G_2 + D < 1$ for $m \leqslant 26$ and
\begin{align*}
G_1 + G_2 &\leqslant \text{Exp}[- (200/101) \pi m \sin(\theta_{0, 3} - \pi / (3 (22 + 6 m))) (\sqrt{3} + 2 \cos(\theta_{0, 3} - 2 \pi / (3 (22 + 6 m))))]\\
 &\quad + \text{Exp}[- (200/51) \pi m \sin(\theta_{0, 3} - \pi / (3 (22 + 6 m))) (\sqrt{3} + 2 \cos(\theta_{0, 3} - 2 \pi / (3 (22 + 6 m))))],
\end{align*}
where the right-hand side is monotonically decreasing in $m$ and is less than $1 - D$ when $m = 27$.

For $k \geqslant 34$, we have $D = 3.8242... \times 10^{-15}$.

When $k \geqslant 34$, for $\theta \in [\theta_{1, 3}, \theta_{0, 3} - \pi / (3 (k + 6 m))]$, we have
\begin{align*}
&\text{Exp}[- 2 \pi m \sin\theta (\sqrt{3} + 2 \cos\theta) / (4 + 2 \sqrt{3} \cos\theta)]\\
 &\leqslant \text{Exp}[- (125/131) \pi m \sin(\pi / (3 (k + 6 m)))]
 \leqslant \text{Exp}[- (333/1048) \pi^2 / ((k / m) + 6)] =: A_1(k/m),\\
&\text{Exp}[- 4 \pi m \sin\theta (\sqrt{3} + 2 \cos\theta) / (7 + 4 \sqrt{3} \cos\theta)]\\
 &\leqslant \text{Exp}[- (250/137) \pi m \sin(\pi t / (k + 6 m))]
 \leqslant \text{Exp}[- (333/548) \pi^2 / ((k / m) + 6)] =: A_2(k/m).
\end{align*}
We also have
\begin{align*}
4 + 2 \sqrt{3} \cos(\theta_{0, 3} - \pi / (3 (k + 6 m))) &\geqslant 1 + (\pi / \sqrt{3}) (1 / (1 + 6 m / k)) (1 / k) =: B_1(k/m),\\
7 + 4 \sqrt{3} \cos(\theta_{0, 3} - \pi / (3 (k + 6 m))) &\geqslant 1 + (2 \pi / \sqrt{3}) (1 / (1 + 6 m / k)) (1 / k) =: B_2(k/m).
\end{align*}

Put $s := k / m$. Then we have $G_1 + G_2 \leqslant A_1(s) B_1(s)^{-k/2} + A_2(s) B_2(s)^{-k/2}$. Here, $A_1$, $A_2$, $B_1$, and $B_2$ are monotonically increasing in $s$. Then,  for $s_1 \leqslant s \leqslant s_2$, we have $G_1 + G_2 \leqslant A_1(s_2) B_1(s_1)^{-k/2} + A_2(s_2) B_2(s_1)^{-k/2}$. The algorithm is as follows, similar to that for the case $p = 2$:

\begin{quotation}
\quad For $s_1 \leqslant s \leqslant s_2$, we determine $a_1$ and $a_2$ so that $B_1(s_1) \geqslant a_1^{2 / k}$ and $B_2(s_1) \geqslant a_2^{2 / k}$, respectively. Then, we have only to show that
\begin{equation*}
(1 / a_1) A_1(s_2) + (1 / a_2) A_2(s_2) < 1 - D.
\end{equation*}
\end{quotation}

We have the following result:
\begin{center}
\begin{tabular}{cc|cc|c}
\hline
$s_1$ & $s_2$ & $a_1$ & $a_2$ & $(1 / a_1) A_1(s_2) + (1 / a_2) A_2(s_2)$\\
\hline
$0$ & $7/20$ & $1$ & $1$ & $0.99914...$\\
$7/20$ & $51/50$ & $21/20$ & $109/100$ & $0.99968...$\\
$51/50$ & $241/100$ & $57/50$ & $31/25$ & $0.99940...$\\
$241/100$ & $551/100$ & $129/100$ & $29/20$ & $0.99989...$\\
$551/100$ & $76/5$ & $153/100$ & $173/100$ & $0.99932...$\\
$76/5$ & $451$ & $189/100$ & $52/25$ & $0.99998...$\\
$451$ & $\infty$ & $239/200$ & $62/25$ & $0.82163...$\\
\hline
\end{tabular}
\end{center}

Finally, when $k \equiv 0 \pmod{12}$, we can write $k = 12 l'$, and we have
\begin{gather*}
Re \, g_{12 l', \; 3}^{*} (\theta_{0, 3} ; R) = e^{(1 / \sqrt{3}) \pi m} \cos((5 l' + m) \pi),\\
Re \, h_{12 l', \; 3}^{*} (\theta_{0, 3} ; R) = e^{(1 / \sqrt{3}) \pi m} \cos((3 l' - m) \pi),\\
Re \, h_{12 l', \; 3}^{* '} (\theta_{0, 3} ; R) = e^{(1 / \sqrt{3}) \pi m} \cos((- l' + m) \pi).
\end{gather*}
Here, the signs of above three terms are in agreement. Furthermore, we have $|F_{12 l', \; 3}^{* ''} (\theta_{0, 3} ; R)| \ll |Re \, g_{12 l', \; 3}^{*} (\theta_{0, 3} ; R)|$.

In conclusion, we can show Lemma \ref{lem-t-m}.\\

\section{Proof of Theorem \ref{th-t+m}}

\subsection{Preliminaries}
Let $k \geqslant 4$ and $m > 0$ be integers, and let $R(t) = t^m$. Then, it is clear that $R$ satisfies the conditions $(i)$, $(ii)$, and $(iii)$ of Property ${P_{k, p}}'$. Furthermore, we have $P_{\gamma_p} = 0$, $N_{\gamma_p} = m$, and $N_R - N_{\gamma_p} = k ((p + 1) / 24) - m$. We may assume that $k ((p + 1) / 24) \geqslant l > m \geqslant 1$.

Similarly to the previous section, we consider a certain number $x$ for the interval $[\theta_{1, p}, \theta_{0, p} - x]$ in which every integer point is included. We have
\begin{equation*}
\arg g_{k, \; p}^{*} (\theta_{0, p} - x ; R) > \arg g_{k, \; p}^{*} (\theta_{0, p} ; R) - (k / 2) x.
\end{equation*}
Thus, when $p = 2$, the bound $x \leqslant \pi / (2 k)$ is sufficient. When $p = 3$ and $k \not\equiv 10 \pmod{12}$, the bound $x \leqslant 2 \pi / (3 k)$ is sufficient.

On the other hand, if $k \not\equiv 10 \pmod{12}$, we need to consider a stricter number $x$ for calculation. We may assume that $x \leqslant \pi / 33$ and $k \geqslant 22$, in which case we have
\begin{equation*}
\arg g_{k, \; 3}^{*} (\theta_{0, 3} - x ; R) > \arg g_{k, \; 3}^{*} (\theta_{0, 3} ; R) - (k / 2 - 9 m / 5) x,
\end{equation*}
and so we may assume $x = \pi / (3 (k - (18 / 5) m))$.

Furthermore, we have  $\left| R \left( e^{2 \pi i (e^{i \theta} / \sqrt{p})} \right) \right| = e^{- (2 / \sqrt{p}) \pi m \sin\theta}$, so we define
\begin{gather}
G := e^{\sqrt{2} \pi m \sin\theta} |h_{k, \; 2}^{*} (\theta ; R)|,\\
G_1 := e^{(2 / \sqrt{3}) \pi m \sin\theta} |h_{k, \; 3}^{*} (\theta ; R)|, \quad
G_2 := e^{(2 / \sqrt{3}) \pi m \sin\theta} |h_{k, \; 3}^{* '} (\theta ; R)|.
\end{gather}

Then, we will show the following lemma in order to prove Theorem \ref{th-t+m}:
\begin{lemma}\label{lem-t+m}\quad
\begin{trivlist}
\item{$(i)$} We have
\begin{gather}
G + (1 / 2) \, {M_{R, 2}}' \, \delta_{k, 2} \, e^{\sqrt{2} \pi m} < 1 \quad \text{for every} \; \theta \in [\theta_{1, 2}, \theta_{0, 2} - \pi / (2 (k + 8 m))], \label{cond-pkp-2} \\
G_1 + G_2 + (1 / 2) \, {M_{R, 3}}' \, \delta_{k, 3} \, e^{(2 / \sqrt{3}) \pi m} < 1 \quad \text{for every} \; \theta \in [\theta_{1, 3}, \theta_{0, 3} - x_0], \label{cond-pkp-3}
\end{gather}
where $x_0 = 2 \pi / (3 k)$ when $k \not\equiv 10 \pmod{12}$ and $x_0 = \pi / (3 (k - (18 / 5) m))$ when $k \equiv 10 \pmod{12}$.

\item{$(ii)$} We have
\begin{equation}
Sign(F_{k, \; p}^{*} (\theta_{0, p} ; R)) = Sign(Re \, g_{k, \; p}^{*} (\theta_{0, p} ; R))
\end{equation}
if $k \equiv 0 \pmod{8}$ for the case of $p = 2$, and if $k \equiv 0 \pmod{12}$ for the case of $p = 3$.
\end{trivlist}
\end{lemma}\quad

\subsection{The case $p = 2$}
We have ${M_{R, 2}}' = 1$ and 
\begin{equation}
(1 / 2) \delta_{k, 2} e^{- \sqrt{2} \pi m}
 \leqslant (1 / 2) e^{- \sqrt{2} \pi} ((1 / 5) e^{\sqrt{2} \pi / 4})^{k / 2} (5 + 10 / (k - 2))
 \leqslant 0.00062185...
\end{equation}
for $k \geqslant 16$ because $(1 / 5) e^{\sqrt{2} \pi / 4} < 1$ and $k \geqslant 8 (m + 1)$.

Furthermore, \begin{equation*}
G = \text{Exp}[2 \pi m \sin\theta (\sqrt{2} + 2 \cos\theta) / (3 + 2 \sqrt{2} \cos\theta)] \; / \; (3 + 2 \sqrt{2} \cos\theta)^{k / 2}.
\end{equation*}
For $\theta \in [\theta_{1, 2}, \theta_{0, 2} - \pi / (2 k)]$, we have
\begin{align*}
1 / (3 + 2 \sqrt{2} \cos(\theta_{0, 2} - \pi / (2 k))) &\leqslant 1 / (1 + \pi / k),\\
\sqrt{2} \sin(\theta_{0, 2} - \pi / (2 k)) &\leqslant 1 + \pi / (2 k),\\
(\sqrt{2} + 2 \cos(\theta_{0, 2} - \pi / (2 k))) / \sqrt{2} &\leqslant (11 / 20) \pi / k.
\end{align*}
We also have
\begin{equation*}
3 + 2 \sqrt{2} \cos(\theta_{0, 2} - \pi / (2 k)) \geqslant 1 + \pi / k \geqslant 4^{2 / k},
\end{equation*}
and so
\begin{equation*}
G \leqslant (1 / 4) e^{(11 / 10) \pi^2 (m / k)} \leqslant (1 / 4) e^{(11 / 80) \pi^2} = 0.97119...
\end{equation*}

Finally, when $k \equiv 0 \pmod{8}$, we can write $k = 8 l$ to yield
\begin{gather*}
Re \, g_{8 l, \; 2}^{*} (\theta_{0, 2} ; R) = e^{- \pi m} \cos((3 l - m) \pi),\\
Re \, h_{8 l, \; 2}^{*} (\theta_{0, 2} ; R) = e^{- \pi m} \cos((l + m) \pi).
\end{gather*}
Here, the signs of above two terms are in agreement. Furthermore, we have $|F_{8 l, \; 2}^{* ''} (\theta_{0, 2} ; R)| \ll |Re \, g_{8 l, \; 2}^{*} (\theta_{0, 2} ; R)|$.

Hence, Lemma \ref{lem-t+m} follows.\\

\subsection{The case $p = 3$}

\subsubsection{Preliminaries}
We have ${M_{R, 3}}' = 1$ and 
\begin{align*}
(1 / 2) \delta_{k, 3} e^{- (2 / \sqrt{3}) \pi m}
 &\leqslant (1 / 2) e^{- (2 / \sqrt{3}) \pi} ((1 / 7) e^{2 \pi / (3 \sqrt{3})})^{k / 2} (13 + 84 / (k - 2))\\
 &\leqslant 0.0034216... \quad =: D
\end{align*}
for $k \geqslant 12$ because we have $(1 / 7) e^{2 \pi / (3 \sqrt{3})} < 1$ and $k \geqslant 6 (m + 1)$.

We have
\begin{align*}
G_1 &= \text{Exp}[2 \pi m \sin\theta (\sqrt{3} + 2 \cos\theta) / (4 + 2 \sqrt{3} \cos\theta)] \; / \; (4 + 2 \sqrt{3} \cos\theta)^{k / 2},\\
G_2 &= \text{Exp}[4 \pi m \sin\theta (\sqrt{3} + 2 \cos\theta) / (7 + 4 \sqrt{3} \cos\theta)] \; / \; (7 + 4 \sqrt{3} \cos\theta)^{k / 2}.
\end{align*}

\subsubsection{The case $k \not\equiv 10 \pmod{12}$}
We may assume that $x = 2 \pi / (3 k)$ and $k \geqslant 12$.

For $12 \leqslant k \leqslant 40$, we have
\begin{equation*}
G_1 + G_2 < \left( (G_1 + G_2) |_{\theta = \theta_{0, 3} - (2 / 3) \pi / k} \right) |_{m = k / 6},
\end{equation*}
and we can show the right-hand side is less than $1 - D$ by numerical calculation.

For $k \geqslant 42$, we have
\begin{align*}
G_1 &\leqslant \text{Exp}[(29/40) \pi^2 m/k] / (1 + 2 \sqrt{3} \pi / (3 k))^{k / 2}\\
G_2 &\leqslant \text{Exp}[(29/20) \pi^2 m/k] / (1 + 4 \sqrt{3} \pi / (3 k))^{k / 2},
\end{align*}
$1 + 2 \sqrt{3} \pi / (3 k) \geqslant (28 / 5)^{2 / k}$, and $1 + 4 \sqrt{3} \pi / (3 k) \geqslant 27^{2 / k}$. Thus, we have
\begin{equation*}
G_1 + G_2 \leqslant (5 / 28) e^{29 \pi^2 / 240} + (1 / 27) e^{29 \pi^2 / 120} = 0.99074...
\end{equation*}

\subsubsection{The case $k \equiv 10 \pmod{12}$}
When $k = 22$, $34$, $46$, and $58$, we have $m \leqslant (k - 10) / 6$, and we can calculate $(G_1 + G_2) |_{\theta = \theta_{0, 3} - \pi / (3 (k - (18 / 5) m))} < 1 - D$ for each $k$ and $m$.

When $k \geqslant 70$, for $\theta \in [\theta_{1, 3}, \theta_{0, 3} - \pi / (3 (k - (18 / 5) m))]$, we have
\begin{align*}
\text{Exp}[2 \pi m \sin\theta (\sqrt{3} + 2 \cos\theta) / (4 + 2 \sqrt{3} \cos\theta)] &\leqslant \text{Exp}[(28/75) \pi^2 / ((k/m) - 18 / 5)] =: A_1(k/m)\\
\text{Exp}[4 \pi m \sin\theta (\sqrt{3} + 2 \cos\theta) / (7 + 4 \sqrt{3} \cos\theta)] &\leqslant \text{Exp}[(56/75) \pi^2 / ((k/m) - 18 / 5)] =: A_2(k/m).
\end{align*}
We also have
\begin{align*}
4 + 2 \sqrt{3} \cos(\theta_{0, 3} - \pi / (3 (k - (18 / 5) m))) &\geqslant 1 + (\pi / \sqrt{3}) (1 / (1 - (18 / 5) m / k)) (1 / k) =: B_1(k/m),\\
7 + 4 \sqrt{3} \cos(\theta_{0, 3} - \pi / (3 (k - (18 / 5) m))) &\geqslant 1 + (2 \pi / \sqrt{3}) (1 / (1 - (18 / 5) m / k)) (1 / k) =: B_2(k/m).
\end{align*}

Following the procedure in the previous section, we put $s := k / m$. $A_1$, $A_2$, $B_1$, and $B_2$ are monotonically decreasing in $s$. Then, for $s_1 \leqslant s \leqslant s_2$, we have $G_1 + G_2 \leqslant A_1(s_1) B_1(s_2)^{-k/2} + A_2(s_1) B_2(s_2)^{-k/2}$. The algorithm is as follows:

\begin{quotation}
\quad For $s_1 \leqslant s \leqslant s_2$, we determine $a_1$ and $a_2$ so that $B_1(s_2) \geqslant a_1^{2 / k}$ and $B_2(s_2) \geqslant a_2^{2 / k}$, respectively. Then, we have only to show that
\begin{equation*}
(1 / a_1) A_1(s_1) + (1 / a_2) A_2(s_1) < 1 - D.
\end{equation*}
\end{quotation}

We have the following result:
\begin{center}
\begin{tabular}{cc|cc|c}
\hline
$s_1$ & $s_2$ & $a_1$ & $a_2$ & $(1 / a_1) A_1(s_1) + (1 / a_2) A_2(s_1)$\\
\hline
$6$ & $25/4$ & $79/10$ & $55$ & $0.97955...$\\
$25/4$ & $20/3$ & $67/10$ & $41$ & $0.99298...$\\
$20/3$ & $22/3$ & $28/5$ & $29$ & $0.97504...$\\
$22/3$ & $35/4$ & $9/2$ & $19$ & $0.97512...$\\
$35/4$ & $13$ & $17/5$ & $11$ & $0.98174...$\\
$13$ & $\infty$ & $12/5$ & $29/5$ & $0.99424...$\\
\hline
\end{tabular}
\end{center}

\subsubsection{The case $k \equiv 0 \pmod{12}$ and $\theta = \theta_{0, 3}$}
We can write $k = 12 l'$, and then we have
\begin{gather*}
Re \, g_{12 l', \; 3}^{*} (\theta_{0, 3} ; R) = e^{- (1 / \sqrt{3}) \pi m} \cos((5 l' - m) \pi),\\
Re \, h_{12 l', \; 3}^{*} (\theta_{0, 3} ; R) = e^{- (1 / \sqrt{3}) \pi m} \cos((3 l' + m) \pi),\\
Re \, h_{12 l', \; 3}^{* '} (\theta_{0, 3} ; R) = e^{- (1 / \sqrt{3}) \pi m} \cos((- l' - m) \pi).
\end{gather*}
Here, the signs of above three terms are in agreement. Furthermore, we have $|F_{12 l', \; 3}^{* ''} (\theta_{0, 3} ; R)| \ll |Re \, g_{12 l', \; 3}^{*} (\theta_{0, 3} ; R)|$.

In conclusion, we have been able to show Lemma \ref{lem-t+m}.\\

\begin{remark}
Note that the location of $m - 1$ zeros is unclear (cf. the end of Section \ref{sec-intro}).

When $l \leqslant 1$, {\it i.e.} ``when $p = 2$, and $4 \leqslant k \leqslant 14, k = 18$'', and ``when $p = 3$, and $4 \leqslant k \leqslant 10, k = 14$'', then we have $m - 1 \leqslant 0$, thus we can prove that all of the zeros of $G_{k, p}^{*}(z ; t^{m})$ in $\mathbb{F}^{*}(p)$ are on the arc $A_p^{*}$ for $p = 2, 3$ and $m \leqslant l$.

When $p = 2$, $k = 16$ and $m = 2$, we can prove that one more zero lies on $A_2^{*}$ by numerical calculation. Thus we can prove that all of the zeros of $G_{16, 2}^{*}(z ; t^{m})$ in $\mathbb{F}^{*}(2)$ lie on $A_2^{*}$ for $m \leqslant 2$. Similarly, when $p = 3$ and $k = 12$, we can also prove that all of the zeros of $G_{12, 3}^{*}(z ; t^{m})$ in $\mathbb{F}^{*}(3)$ lie on $A_3^{*}$ for $m \leqslant 2$ by numerical calculation. 
\end{remark}\quad

\begin{center}
{\large\sc Acknowledgement.}
\end{center}
The author thanks Professor Masanobu Kaneko and Professor Eiichi Bannai for suggesting these problems.\\

\quad\\

\begin{flushright}\makebox{
\begin{minipage}{2.8in}
Faculty of Mathematics\\
Graduate School\\
Kyushu University\\
Motooka 744, Nishi-ku\\
Fukuoka, 819-0395, Japan\\
j.shigezumi@math.kyushu-u.ac.jp
\end{minipage}}\end{flushright}

\begin{thebibliography}{9}

\bibitem{AKN}
T. Asai, M. Kaneko, and H. Ninomiya: {\it Zeros of certain modular functions and an application}, Comment.
Math. Univ. St. Paul. {\bfseries 46} (1997), 93--101.

\bibitem{BKM}
E. Bannai, K. Kojima, and T. Miezaki: {\it On the zeros of Hecke type Faber polynomial}, Kyushu. J. Math. {\bfseries 62}(2007), 15--61.

\bibitem{DJ}
W. Duke and P. Jenkins: {\it On the zeros and coefficients of certain weakly holomorphic modular forms}, Pure Appl. Math. Q. {\bfseries 4} (2008), No. 4, 1327--1340.

\bibitem{Ka}
M. Kaneko, {\it On the zeros of certain modular forms}, in: Number theory and its applications (Kyoto, 1997), 193--197, Dev. Math., 2, Kluwer Acad. Publ., Dordrecht, 1999.

\bibitem{Kr}
A. Krieg: {\it Modular Forms on the Fricke Group.}, Abh. Math. Sem. Univ. Hamburg, {\bfseries 65}(1995), 293-299.

\bibitem{MNS}
T. Miezaki, H. Nozaki, and J. Shigezumi: {\it On the zeros of Eisenstein series for $\Gamma_0^{*}(2)$ and $\Gamma_0^{*}(3)$}, J. Math. Soc. Japan, {\bfseries 59}(2007), 693--706.

\bibitem{Q}
H. -G. Quebbemann: {\it Atkin-Lehner eigenforms and strongly modular lattices.}, Enseign. Math. (2), {\bfseries 43}(1997), No. 1-2, 55-65.

\bibitem{R}
R. A. Rankin: {\it The zeros of certain Poincar\'{e} series}, Compositio Math., {\bfseries 46}(1982), 255-272.

\bibitem{RSD}
F. K. C. Rankin and H. P. F. Swinnerton-Dyer: {\it On the zeros of Eisenstein Series}, Bull. London Math. Soc., {\bfseries 2}(1970), 169-170.

\bibitem{T}
H. Tokitsu: {\it On the location of the zeros of certain modular forms $($Japanese$)$}, M. S. thesis, (2008), Kyushu University.

\end{thebibliography}
\end{document}